\input amstex
\documentstyle{amsppt}
\topmatter
\title Projections in normed linear spaces and
 sufficient enlargements
\endtitle
\rightheadtext{sufficient enlargements}
\author M.I.Ostrovskii
\endauthor

\address  Department of Mathematics,
The Catholic University of America, 
Washington, D.C. 20064, USA
\endaddress

\email ostrovskii\@cua.edu\endemail

\abstract 
{\bf Definition.} A symmetric with respect to $0$ bounded closed convex 
set $A$ in a finite dimensional normed space $X$ is called a {\it 
sufficient enlargement} for $X$ (or of $B(X)$) if for arbitrary
isometric embedding of $X$ into a Banach space $Y$ there exists a projection
$P:Y\to X$ such that $P(B(Y))\subset A$ (by $B$ we denote the unit ball). 

The notion of sufficient enlargement is implicit in the paper: B.Gr\"unbaum, 
Projection constants, Trans. Amer. Math. Soc. {\bf 95} (1960) 451--465. 
It was explicilty introduced by the author in: M.I.Ostrovskii,
Generalization of projection constants: sufficient enlargements, Extracta 
Math., {\bf 11} (1996), 466--474.

The main purpose of the present paper is to continue investigation of 
sufficient enlargements started in the papers cited above. In particular the 
author investigate sufficient enlargements whose support functions are in 
some directions close to those of the unit ball of the  space, sufficient 
enlargements of minimal volume,  sufficient enlargements for euclidean spaces.
\endabstract
\thanks The research was supported by an INTAS
grant and by a grant of T\"UBITAK, 
the first version of the paper was prepared
when the author was visiting the University of Michigan
(Ann Arbor) and Odense University. 
The author would like to thank N.J.Nielsen and M.S.Ramanujan
for their hospitality. The author is obliged to the referees
for many useful remarks.\endthanks

\endtopmatter

\head \S 1. Introduction \endhead

We denote the unit ball (sphere) of a normed linear space $X$ by 
$B(X)$ ($S(X)$).

{\bf Convention.} We shall use the term {\it ball} \ for 
symmetric with respect to $0$ 
bounded closed convex set with nonempty interior
in a finite dimensional linear space.

\definition{Definition 1} A ball $A$ in a finite dimensional normed
space $X$ is called a {\it sufficient enlargement} for $X$ (or of $B(X)$)
if for arbitrary
isometric embedding $X\subset Y$ ($Y$ is a Banach space)
there exists a projection
$P:Y\to X$ such that $P(B(Y))\subset A$. A {\it minimal sufficient
enlargement} is defined to be a sufficient enlargement no proper
subset of which is a sufficient enlargement.
\enddefinition

The notion of sufficient enlargement is 
implicit in B.Gr\"unbaum's paper \cite{2}, it was explicilty
introduced by the present author in \cite{5}.

The notion of sufficient enlargement is of interest because
it is a natural geometric notion, it characterizes
possible shadows of symmetric convex body onto
a subspace, whose intersection with the body is given.

The main purpose of the present paper is to continue
investigation of sufficient enlargements started in \cite{5}.
In \S 2 we investigate sufficient enlargements whose support
functions are in some directions close to those of the
unit ball of the space, \S 3 is devoted to sufficient
enlargements for euclidean spaces.

We refer to \cite{4} and \cite{7} for background on
Banach space theory and to \cite{6} for background on the
theory of convex bodies.

{\bf Some recalls.}
Let $X$ and $Y$ be finite dimensional normed spaces
and $T:X\to Y$ be a linear operator. An {\it $l_\infty-$factorization} of $T$
is a pair of operators $u_1:X\to l_\infty$ and $u_2:l_\infty\to Y$
satisfying $T=u_2u_1$. The {\it $L_\infty-$factorable
norm} of $T$ is defined to be the $\inf||u_1||||u_2||$,
where the $\inf$ is taken over all $l_\infty-$factorizations.

An {\it absolute projection constant} of a finite dimensional normed linear
space $X$ is defined to be the smallest positive real number
$\lambda(X)$ such that for every isometric embedding $X\subset Y$
there exists a continuous linear projection $P:Y\to X$ with $||P||\le 
\lambda(X)$.

We shall use the following observations.

\proclaim{Proposition 1} \cite{5}. Let $A$ be a  ball in a finite dimensional
normed linear space $X$.
The space $X$ normed by the gauge functional
of $A$ will be denoted by $X_A$. 

The ball $A$ is a sufficient enlargement
for $X$ if and only if the $L_\infty-$factorable norm of
the natural identity mapping from $X$ to $X_A$ is $\le 1$.
\endproclaim

\proclaim{Proposition 2} \cite{2}. A symmetric with respect to $0$
parallelepiped containing $B(X)$
is a sufficient enlargement for $X$.
\endproclaim

\proclaim{Proposition 3} \cite{2}. Convex combination of 
sufficient enlargements 
for $X$ is a sufficient
enlargement for $X$. 
\endproclaim

\head \S 2. Sufficient enlargements which are in some directions
close to the balls\endhead

We start with an investigation of a sufficient enlargement
which is contained in a homothetic image of a circumscribed parallelepiped
with the coefficient of homothety close to 1 (and of course
greater than 1). Next result gives a condition under
which such enlargement  contains a non-trivial
homothetic image of the parallelepiped.
\bigskip

\proclaim{Theorem 1} Let $X$ be an $n-$dimensional normed
space. Let $\{f_i\}_{i=1}^n\subset S(X^*)$ be a basis of $X^*$
and let vectors
$x_i\in S(X)$ be such that $f_i(x_i)=1$ and for
some $c_2>0$ and each $f\in B(X^*)$ there exists  at most one
element $i$ in the set $\{1,\dots,n\}$ for which
$|f(x_i)|\ge 1-c_2.$

Let $A$ be a sufficient enlargement for $X$ such
that for some $c_1\ge 0$ it is
contained in the parallelepiped
$\{x:\ |f_i(x)|\le1+c_1,\ i\in\{1,\dots,n\}\}$

Let
$c_3=1-\frac{2-c_2}{c_2}c_1.$ Suppose $c_3>0$.
Then $A$ contains the parallelepiped
$Q:=\{x:\ |f_i(x)|\le c_3,\ i\in\{1,\dots,n\}\}.$
\endproclaim

\demo{Proof} Let $\{f_i\}_{i=n+1}^\infty\subset S(X^*)$ be such that
$(\forall x\in X)\ (||x||=\sup\{|f_i(x)|:\ i\in\Bbb N\}).$
Then the operator
$E:X\to l_\infty$ 
defined by
$Ex:=\{f_i(x)\}_{i=1}^\infty$ is an isometric embedding.
Let $P:l_\infty\to E(X)$ be a projection for
which $P(B(l_\infty))\subset E(A).$

The condition of the theorem imply that
there exists a partition of $\Bbb N$ into
subsets $F_1,\dots, F_n$ such that for $i\in F_j$
we have $f_i(x_k)<1-c_2$ for $k\ne j$. 

Let us show that $P(B(l_\infty))$ contains $E(Q)$.
Observe that the first $n$ coordinate functionals
on $l_\infty$ are norm-preserving extensions of
functionals $f_iE^{-1}:E(X)\to\Bbb R$. Therefore
in order to prove that $A\supset Q$ it is sufficient
to prove that for every collection $\{\theta_i\}_{i=1}^n,\
\theta_1=\pm 1$ there exists a vector $z_\theta\in 
B(l_\infty)$ and real numbers $b_1,\dots,b_n\ge c_3$
such that
$$Pz_\theta=(\theta_1b_1,\theta_2b_2,\dots,
\theta_nb_n, b_{n+1}, b_{n+2},\dots)$$
for some $b_{n+1}, b_{n+2},\dots\in\Bbb R$.

We introduce $z_\theta$ as the sequence
$\{d_k\}_{k=1}^\infty$, where $d_k=\theta_jf_k(x_j)$
if $k\in F_j$. In particular, $d_1=\theta_1,\dots,
d_n=\theta_n$. Let us show that $Pz_\theta$
satisfies the requirement above. Let
$$Pz_\theta=(\alpha_1,\dots,\alpha_n,\alpha_{n+1},\dots).$$
Suppose that for some $m\in\{1,\dots,n\}$ we have
$\alpha_m\notin[\theta_mc_3,\theta_m\infty)$.
Let us consider the family of vectors
$$y_\delta=(1+\delta)\theta_mE(x_m)-
\delta z_\theta,\ (\delta>0).$$

When $\delta>0$ is small enough, then $y_\delta
\in B(l_\infty)$. More precisely, by the conditions
of the theorem it happens at least when
$(1-c_2)(1+\delta)+\delta\le 1,$ that is, when
$\delta\le\frac{c_2}{2-c_2}.$

On the other hand the $m-$th coordinate of
$Py_\delta$ is equal to
$$(1+\delta)\theta_m-\delta\alpha_m=
\theta_m+\delta(\theta_m-\alpha_m).$$
So for $0\le\delta\le c_2/(2-c_2)$ we have
$|\theta_m+\delta(\theta_m-\alpha_m)|\le 1+c_1.$
Hence
$$1+\frac{c_2}{2-c_2}(1-c_3)<1+c_1\
\hbox{ or }\ 
c_3>1-\frac{2-c_2}{c_2}c_1.$$
This contradicts the condition on $c_3$.
$\ \ \square$
\enddemo
\bigskip

\proclaim{Corollary} Let $X$ be an $n$-dimensional normed space
 and $Q$ be a parallelepiped circumscribed about $B(X)$. 
Suppose there exist points
 $\{x_i\}_{i=1}^n$ on faces of $Q$ (one point on the union
 of each pair of symmetric faces) such that  $x_i\in B(X)$
 and for every pair $(x_i,x_j),\ x_i\ne x_j$ and every
 $f\in B(X^*)$ at least one of the numbers $|f(x_i)|$ is 
less than $1$. Then
$Q$ is a minimal sufficient enlargement for $X$. 
\endproclaim
 
\demo{Proof} By Proposition 2 only minimality requires a
proof. Let $\{f_i\}_{i=1}^n\subset B(X^*)$ be such that
$Q=\{x:\ |f_i(x)|\le 1,\ i\in\{1,\dots,n\}\}$.

By compactness of $B(X^*)$ there exists $c_2>0$ satisfying the
condition of Theorem 1. Let $A\subset Q$ be a sufficient enlargement
for $X$. Applying Theorem 1 with $c_1=0$ we get $A\supset Q$.
Hence the sufficient enlargement $Q$ is minimal. $\square$
\enddemo

\remark{Remark 1} Condition $|f(x_i)|\ge 1-c_2$ in Theorem 1
cannot be omitted. This statement can be derived e.g. from
the following observation which is interesting itself:
the proof of the M.Kadets--Snobar theorem as it is given
in \cite{3}, (see, also \cite{7}, \S 15) shows the following.
Let $X$ be an $n-$dimensional normed linear space and
$E\subset B(X)$ be the ellipsoid of maximal volume in
$B(X)$. Then $\sqrt{n}E$ is a sufficient enlargement for
$X$. In particular $B(l_2^n)$ is a sufficient enlargement
for $B(l_1^n).$ Letting $A=B(l_\infty^n)$ we get the
statement.
\endremark

\remark{Remark 2} The following example shows that there are
no direct generalizations of Theorem 1 for non-trivially
large values of $c_1$:

For arbitrary $h\in S(l_2^n)$ there exists a sufficient
enlargement for $l_2^n$ which is contained in the intersection
of $3B(l_\infty^n)$ and the set $\{x:\ |\langle h,x\rangle|\le 
1\}$. 

In fact, let $P_h$ be a projection onto the
hyperplane orthogonal to $h$ with minimal possible norm as
an operator on
$l_\infty^n$ and let
$A=[-h,h]+P_h(B(l_\infty^n)).$
It is easy to see that $A$ is a sufficient enlargement for
$l_2^n$ satisfying all the requirements.
\endremark

\bigskip

The next result shows that the condition of the Corollary is
not necessary for $Q$ to be a minimal sufficient
enlargement.

\proclaim{Theorem 2} 
There exist a two-dimensional normed linear space $X$ and 
functionals $f_1, f_2\in B(X^*)$
such that the following conditions are satisfied:

1) There exists precisely one point $x_1\in B(X)$ such that
$f_1(x_1)=1$ and precisely one point $x_2\in B(X)$ such that
$f_2(x_2)=1$.

2) The parallelogram 
$C=\{x:\ |f_1(x)|\le 1,\ |f_2(x)|\le 1\}$
is a minimal sufficient enlargement.

3) There exist a linear functional $f_3\in B(X^*)$ such
that $|f_3(x_1)|=|f_3(x_2)|=1$.
\endproclaim

\demo{Proof} Consider the space whose unit ball is the euclidean
disc intersected with the strip
$$\{(a_1,a_2):\ |a_1-a_2|\le 1\}.$$

Let $x_1=(1,0),\ x_2=(0,1)$ and let
$f_1$ and $f_2$
be the coordinate functionals. It is clear that 
Condition 1 of the theorem is satisfied.

In our case
$C=\{(a_1,a_2):\ |a_1|\le 1,\ |a_2|\le 1\}.$

It is clear that the functional 
$f_3(a_1,a_2)=a_1-a_2$ satisfies Condition 3 of the
theorem.

It remains to show,
that  $C$ is a minimal sufficient enlargement.

Let $\{f_i\}_{i=4}^\infty\subset S(X^*)$ be such that
$(\forall x\in X)\ (||x||=\sup\{|f_i(x)|:\ i\in\Bbb N\}).$
Then the operator
$E:X\to l_\infty$ 
defined by
$Ex:=\{f_i(x)\}_{i=1}^\infty$ is an isometric embedding.

Now, if we suppose that $C$ is not a minimal sufficient
enlargement, then there exists a projection $P:l_\infty\to
E(X)$, such that the closure of
its image is a proper part of $E(C)$.
We show that this gives us a contradiction.

Consider the vectors
$$x_1(\varepsilon):=(\cos\varepsilon, \sin\varepsilon),\
x_2(\varepsilon):=(\sin\varepsilon, \cos\varepsilon)\in B(X),\
0<\varepsilon<\pi/4.$$

It is clear that for $0<\varepsilon<\pi/4$
the following is true (the reader is advised to 
draw the picture): 
for each $f\in B(X^*)$ either
$$|f(x_1(\varepsilon))|\le 1-\tan\varepsilon\
\hbox{ or }\
|f(x_2(\varepsilon))|\le 1-\tan\varepsilon.$$

Therefore there exists a partition 
$\Bbb N=A_1(\varepsilon)\cup A_2(\varepsilon)$
such that
$|f_i(x_1(\varepsilon))|\le 1-\tan\varepsilon$ 
for $i\in A_2(\varepsilon)$  and
$|f_i(x_2(\varepsilon))|\le 1-\tan\varepsilon$
for $i\in A_1(\varepsilon)$.

Now for $\theta=(\theta_1,\theta_2),$ where $\theta_1=\pm 1,
\ \theta_2=\pm 1,$ we define $z_\theta(\varepsilon)\in
l_\infty$ as the vector, whose $i-$th coordinates coincide
with the coordinates of $\theta_1 Ex_1(\varepsilon)$ for
$i\in A_1(\varepsilon)$ and with the coordinates of
$\theta_2E x_2(\varepsilon)$ for $i\in A_2(\varepsilon)$.

It is clear that $z\in B(l_\infty)$. Let
$$Pz_\theta(\varepsilon)=(\alpha_1,\alpha_2,\dots,
\alpha_n, \dots)\in l_\infty.$$

Let us show that

$$\theta_1\alpha_1\ge\cos\varepsilon-2(1-\cos\varepsilon)/
\varepsilon,\eqno{(1)}$$
$$\theta_2\alpha_2\ge\cos\varepsilon-2(1-\cos\varepsilon)/
\varepsilon.\eqno{(2)}$$

Because $\varepsilon>0$ and
$\theta=(\theta_1,\theta_2)=(\pm 1,\pm 1)$ are arbitrary
(1) and (2) imply
$P(B(l_\infty))\supset E(C),$
so we get a contradiction.

Suppose that either (1) or (2) is not satisfied.
Without loss of generality, we may assume that (1)
is not satisfied.

Consider the family of vectors

$$y_\delta=(1+\delta)\theta_1E(x_1(\varepsilon))-
\delta z_\theta(\varepsilon)\in l_\infty\ 
(\delta>0).$$

From the definition of $z_\theta(\varepsilon)$ it
is easy to derive that
$$||y_\delta||_\infty\le\max\{1,(1+\delta)
(1-\tan\varepsilon)+\delta\}.$$

Hence if $\delta$ is such that $2\delta/(1+\delta)\le
\tan\varepsilon$, then $||y_\delta||_\infty\le 1$.
In particular, $||y_{\varepsilon/2}||_\infty\le 1$.
Since $P(B(l_\infty))\subset E(C)$, then
the modulus of the first coordinate of 
$Py_{\varepsilon/2}\in
l_\infty$ is $\le 1$.
On the other hand, we have

$$Py_{\varepsilon/2}=
(1+\varepsilon/2)\theta_1E(x_1(\varepsilon))-
(\varepsilon/2)Pz_\theta(\varepsilon).$$
Hence the first coordinate of $Py_{\varepsilon/2}$
is
$$(1+\varepsilon/2)\theta_1\cos\varepsilon-
(\varepsilon/2)\alpha_1.$$
We have
$$|(1+\varepsilon/2)\theta_1\cos\varepsilon-
(\varepsilon/2)\alpha_1|=
|(1+\varepsilon/2)\cos\varepsilon-
(\varepsilon/2)\theta_1\alpha_1|>$$
$$(1+\varepsilon/2)\cos\varepsilon-
(\varepsilon/2)(\cos\varepsilon-
2(1-\cos\varepsilon)/\varepsilon)=1.$$

This contradiction implies that (1) and (2) are valid.
Theorem 2 is proved. $\square$
\enddemo

By a {\it prism} in $\Bbb R^n$ we mean the Minkowski sum of a
set $A$ lying in an $(n-1)-$dimensional hyperplane and a line segment
that is not parallel to the hyperplane. The set $A$ is called a {\it basis}
of the prism.
\smallskip

It turns out that if a sufficient enlargement $A$ for $X$
is such that its boundary intersects
$S(X)$ in a smooth point, then $A$ should contain a prism,
which is also a sufficient enlargement, so the investigation
of such enlargement can be in certain sense reduced to investigation
of $(n-1)-$dimensional sufficient enlargement.

\proclaim{Theorem 3} Let $X$ be an $n-$dimensional normed space
and let $x_1\in S(X)$ be a smooth point and $h\in S(X^*)$ be
its supporting functional. Let $\{x_i\}_{i=2}^n\subset S(X)$
be such that $\{x_i\}_{i=1}^n$ is a basis in $X$ and
$h(x_i)=0$ for $i\in\{2,\dots,n\}$. Suppose that
$A$ is a sufficient enlargement
for $X$, which is contained in the set
$\{x\in X:\ |h(x)|\le 1\}.$
Then there exists a symmetric with respect to $0$ prism
$M$ with basis parallel to {\rm $\hbox{lin}\{x_2,\dots,x_n\}$
}such that

(a) $M\subset A$;

(b) $M$ is a sufficient enlargement for $X$.
\endproclaim

\demo{Proof} We consider the natural isometric
embedding $E$ of $X$ into $C(S(X^*))$: every vector
is mapped onto its restriction  (as a function
on $X^*$) to $S(X^*)$. We introduce the following notation:
$C=C(S(X^*))$ and $B_C=B(C(S(X^*)))$.

Since $A$ is a sufficient enlargement for
$X$, then there exists a projection
$P: C\to \hbox{lin}\{Ex_i\}_{i=1}^n,$
such that
$$P(B_C)\subset E(A).
\eqno{(3)}$$

Projection $P$ can be represented as
$P(f)=\sum_{i=1}^n\mu_i(f)Ex_i,$
where $\mu_i$ are measures on $S(X^*)$.

Inclusion (3) implies that $||\mu_1||\le 1$. Since
$P$ is a projection we have
$\mu_j(Ex_i)=\delta_{i,j}\ (i,j=1,\dots,n).$
In particular, $\mu_1(Ex_1)=1$. Because $x_1$ is
a smooth point, the function $|Ex_1|\in C$ attains its maximum
only at $h$ and $-h$. Hence   $\mu_1$ can be represented as
$\mu_1=b_{1,1}\delta_{h}+b_{2,1}\delta_{-h},$
where $\delta_h$ and $\delta_{-h}$ are Dirac measures,
$b_{1,1}\ge 0,\ b_{2,1}\le 0$ and $b_{1,1}-
b_{2,1}=1$.

Now, for $i=2,\dots,n$ we find representations
$$\mu_i=b_{1,i}\delta_{h}+b_{2,i}\delta_{-h}+
\nu_i,$$
where $\nu_i$ don't have atoms in $h$ and $-h$.
To unify the notation we set $\nu_1=0$.

We introduce new measures
$$\omega_i:=(b_{1,i}-b_{2,i})\delta_{h}+\nu_i.$$
It is clear that
$\omega_j(Ex_i)=\delta_{i,j}\ (i,j=1,\dots,n).$
Hence
$Q(f):=\sum_{i=1}^n \omega_i(f)Ex_i$
is also a projection onto lin$\{Ex_i\}_{i=1}^n$.

Let us show that
$$Q(B_C)\subset\hbox{cl}(P(B_C)).
\eqno{(4)}$$

Let $f\in B_C$.
Since $\nu_i$ don't have atoms in $\pm h$, then for
every $\varepsilon>0$ there exists a function
$g\in B_C$ such that $g(-h)=-f(h)$,
$g(h)=f(h)$ and $|\nu_i(f)-\nu_i(g)|<\varepsilon$
for all $i\in\{1,\dots,n\}$. This implies that
$$\forall i\in\{1,\dots,n\}\ |\omega_i(f)-
\mu_i(g)|<\varepsilon.$$
Since $\varepsilon>0$ is arbitrary (4) follows.
Hence
$Q(B_C)\subset E(A).$
Now we shall show that
$M:=E^{-1}(\hbox{cl}(Q(B_C)))$
is the required prism.

The fact that $M$ is a sufficient enlargement
follows by a standard argument from the fact that $C$ is an
$\Cal L_\infty$-space (see \cite{4}).

It remains to show that $E(M)$ is a prism with 
basis parallel to lin$\{Ex_2,\dots,Ex_n\}$.

We have
$$E(M)=\hbox{cl}\{f(h)Ex_1+
\sum_{i=2}^n(b_{1,i}-b_{2,i})f(h)Ex_i+
\sum_{i=2}^n\nu_i(f)Ex_i:\ f\in B_C\}.$$

It is clear that the closures of the sets
$$\Gamma_\alpha:=\{\sum_{i=2}^n\nu_i(f)Ex_i:
\ f\in B_C,\ f(h)=\alpha\}$$
don't depend on $\alpha$. So $M$ is a prism
of required form.
The theorem is proved. $\square$ 
\enddemo

\head \S 3. Sufficient enlargements for euclidean spaces
\endhead

Dealing with sufficient enlargements for $l_2^n$ it 
is useful to introduce the following definition.

\definition{Definition 2} A sufficient enlargement
$A$ for $l_2^n$ is said to be {\it small} if
$$\int_{O(n)}T(A)d\mu(T)=
\lambda(l_2^n)B(l_2^n),$$
where $\mu$ is the normalized Haar measure on the
orthogonal group $O(n)$ and $\lambda(l_2^n)$ is the absolute
projection constant.
\enddefinition

\remark{Remark} Proposition 3 implies that
$$\int_{O(n)}T(A)d\mu(T)\supset
\lambda(l_2^n)B(l_2^n)$$
for arbitrary sufficient enlargement $A$.
This explains the choice of the term ``small''.
\endremark

The following result supplies us with a wide and interesting
class of small sufficient enlargements.

\proclaim{Theorem 4} Let $G$ be a finite subgroup of $O(n)$
such that each linear operator on $\Bbb R^n$ commuting with all
elements of $G$ is a scalar multiple of the identity.
Then for every $y\in S(l_2^n)$ the Minkowski sum of segments
$$A=\frac{n}{|G|}\sum_{g\in G}[-g(y),g(y)]$$
is a small sufficient enlargement for $l_2^n$.
\endproclaim

\demo{Proof} First we prove
$$(\forall x\in\Bbb R^n)\ (x=\frac{n}{|G|}\sum_{g\in G}
\langle x, g(y)\rangle g(y)).\eqno{(5)}$$
Let us introduce a linear operator
$T:l_2^n\to l_2^n$ by the equality
$$Tx=\sum_{g\in G}\langle x, g(y)\rangle g(y)
\eqno{(6)}$$

Let us show that $hT=Th$ for each $h\in G$. In fact
$$hT(x)=\sum_{g\in G}\langle x, g(y)\rangle hg(y)=
\sum_{g\in G}\langle h(x), hg(y)\rangle hg(y)=$$
$$\sum_{g\in G}\langle h(x), g(y)\rangle g(y)=Th(x).$$
Hence $T=\lambda I$ for some
$\lambda\in\Bbb R$.

The equality of traces in (6) shows that $\lambda n=|G|$.
Hence $\lambda=\frac{|G|}n$. The assertion (5) follows.

Now, (5)
implies that the identity operator on $l_2^n$ admits
factorization
$I=T_2T_1,$
where
$T_1:l_2^n\to l_\infty^G$ and  $T_2:l_\infty^G\to l_2^n$
are defined as follows
$$T_1(x)=\{\langle x,g(y)\rangle\}_{g\in G}\
\hbox{ and }\
T_2(\{a_g\}_{g\in G})=\frac{n}{|G|}\sum_{g\in G}a_g g(y).$$
It is clear that $||T_1||=1$ and $A=T_2(B(l_\infty^G))$, 
therefore $A$ is a sufficient enlargement (see Proposition
1).
\bigskip

The enlargement $A$ is small by
the following observation. A calculation of B.Gr\"unbaum
\cite{2} shows that 
$$\forall z\in l_2^n\
\int_{O(n)}T([-z,z])d\mu(T)=\frac{||z||\lambda(l_2^n)}nB(l_2^n).
\eqno{(7)}$$
Therefore
$$\int_{O(n)}T(A)d\mu(T)=\frac{n}{|G|}\sum_{g\in G}
\frac{||g(y)||\lambda(l_2^n)}nB(l_2^n)=
\lambda(l_2^n)B(l_2^n).$$
$\square$
\enddemo

\remark{Remark 1}
It is easy to find examples showing that for different $y\in S(l_2^n)$
we get quite different sufficient enlargements.
\endremark

\remark{Remark 2} Many different groups satisfying the condition of
Theorem 4 are given by the representation theory of finite groups.
In particular, every irreducible real representation 
of a finite group
satisfies the condition (after a proper choice of an inner product on
$\Bbb R^n$).
\endremark
\medskip

Small sufficient enlargements have the following nice property.

\proclaim{Theorem 5} Let $A$ be a sufficient enlargement
for $l_2^{n+m}=l_2^n\oplus l_2^m$ and suppose that the images
$A_1$ and $A_2$ of $A$ by
the orthogonal projections onto $l_2^n$ and $l_2^m$ are
small sufficient enlargements for $l_2^n$ and 
$l_2^m$.
Then $A=A_1+A_2$ (Minkowski sum).
\endproclaim

\demo{Proof} We claim: if $A_1$
and $A_2$ are small sufficient enlargements for
$l_2^n$ and $l_2^m$, then $A_1+A_2\subset l_2^{n+m}$ is
a small sufficient enlargement.

At the moment we do not need the fact that $A_1+A_2$ is a
sufficient enlargement, but because the proof is simple,
we sketch it. By Proposition 1 the fact that
$A_1$ is a sufficient enlargement for $l_2^n$ means that
the $L_\infty-$factorable norm of the
identical embedding of $l_2^n$ into $\Bbb R^n$ normed by
the gauge functional of $A_1$ is not greater than $1$, the
analogous assertion is valid for $l_2^m$ and $A_2$. Now,
it is easy to see that the $L_\infty-$factorable norm of
the identical embedding of $l_2^n\oplus_2l_2^m$ into
$\Bbb R^{n+m}$ normed by the gauge functional of
$A_1+A_2$ is $\le 1$.
\bigskip

The fact that the sufficient enlargement $A_1+A_2$ is small
can be proved in the following way:
$$\int_{O(n+m)}T(A_1+A_2)d\mu(T)=$$
$$\int_{O(n+m)}T(\int_{O(n)}T_1(A_1)d\mu_1(T_1)+
\int_{O(m)}T_2(A_2)d\mu_2(T_2))d\mu(T)=$$
(here $\mu_1$ and $\mu_2$ are normalized Haar measures on
$O(n)$ and $O(m)$ respectively)
$$\int_{O(n+m)}T(\lambda(l_2^n)B(l_2^n)+
\lambda(l_2^m)B(l_2^m))d\mu(T)=$$
$$\int_{O(n+m)}T(\int_{O(n)}T_1(Q_1)d\mu_1(T_1)+
\int_{O(m)}T_2(Q_2)d\mu_2(T_2))d\mu(T)=$$
(here $Q_1$ and $Q_2$ are cubes cicumscribed about $B(l_2^n)$
and $B(l_2^m)$ respectively)
$$\int_{O(n+m)}T(Q_1+Q_2)d\mu(T)=
\lambda(l_2^{n+m})B(l_2^{n+m})$$
(by B.Gr\"unbaum's result \cite{2}).

Now the theorem follows from the following direct consequence
of the remark after Definition 2: small sufficient
enlargements are minimal, in particular, no proper 
subset of $A_1+A_2$ is a sufficient enlargement
for $l_2^{n+m}$.
$\square$
\enddemo

Let $X$ be a finite dimensional normed linear space.
Denote by {\bf M} the set of all sufficient enlargements
of minimal volume for $X$. Results of \cite{1} (Theorem 6)
imply the following result.

\proclaim{Theorem 6} The set {\bf M} contains a parallelepiped.
\endproclaim

Easy examples (e.g. two-dimensional space whose ball is 
regular hexagon) show that {\bf M} may contain balls which
are not parallelepipeds. But
it turns out that for Euclidean spaces {\bf M} contains only
parallelepipeds.

\proclaim{Theorem 7} If $A$ is a sufficient enlargement of
minimal volume for $l_2^n$, then $A$ is a cube
circumscribed about $B(l_2^n)$.
\endproclaim

\demo{Proof} Let $A$ be a sufficient
enlargement for $l_2^n$ and vol$A=2^n$. We may assume
without loss of generality (see Proposition 1)
that $A$ is a zonoid. Therefore (see \cite{6}, p. 183), its
support function can be represented in the form
$$h(A,x)=\int_{S^{n-1}}|\langle x,v\rangle|d\rho(v)
\hbox{ for } x\in\Bbb R^n$$
with some even measure $\rho$ on $S^{n-1}$.

We denote by $D$ the set of all smooth points on the
boundary of $A$. It is known (see \cite{6}, p. 73) that
the complement of $D$ in the surface of $A$ has zero
surface measure. Let $T:D\to S^{n-1}$ be the spherical
image map (see \cite{6}, p. 78), that is: 
$T(d)$ is the unique outer unit normal vector
of $A$ at $d$. Let $\mu$ be the measure on $S^{n-1}$
defined by
$$\mu(\Omega)=m_{n-1}(T^{-1}(\Omega)),$$
where $m_{n-1}$ is the surface area measure on the boundary
of $A$.

It is clear that
$$\hbox{vol}A=\frac1n\int_{S^{n-1}}h(A,x)d\mu(x)=
\frac1n\int_{S^{n-1}}\int_{S^{n-1}}
|\langle x,v\rangle|d\rho(v)d\mu(x).$$

The $(n-1)-$dimensional volume of the orthogonal
projection of $A$ onto the hyperplane orthogonal
to $w\in S^{n-1}$ can be computed as
$$\alpha(w)=\frac12\int_{S^{n-1}}|\langle x,w
\rangle|d\mu(x).$$

We proceed by induction on the dimension. The case
$n=1$ is trivial. Suppose that we have  proved the result for
$n-1$. Now, let $A$ be a sufficient enlargement for
$l_2^n$ and vol$A=2^n$. 

By Fubini theorem
$$2^n=\hbox{vol}A=\frac1n\int_{S^{n-1}}2\alpha(w)d\rho(w).$$

Since $A$ is a sufficient enlargement, it is easy to derive
from (7)
that var$(\rho)\ge n$.

It is clear that an orthogonal projection of $A$ onto an 
$(n-1)-$dimensional subspace is a sufficient enlargement
for $l_2^{n-1}$. It is clear also that every parallelepiped
containing $B(l_2^{n-1})$ has volume $\ge 2^{n-1}$. Therefore
by Theorem 6 $\alpha(w)\ge 2^{n-1}$. 
It follows that almost everywhere (in the sense of $\rho$)
$\alpha(w)=2^{n-1}$.

By induction hypothesis orthogonal
projections in directions $w$ for which 
$\alpha(w)=2^{n-1}$ are cubes. Let us choose
one such direction, say $w_1$, and let us denote by
$w_2,w_3,\dots,w_n$ an orthonormal basis in the
subspace orthogonal to $w_1$ such that the orthogonal
projection of $A$ onto lin$\{w_2,\dots,w_n\}$ is
$$[-w_2,w_2]+\dots+[-w_n,w_n].$$
In particular
$$A\subset\{x:|\langle x,w_2\rangle|\le 1\}.$$
By Theorem 3 $A$ contains a prism $M$ with the basis
parallel to
$$\hbox{lin}\{w_1, w_3,w_4,\dots,w_n\}$$
such that
$M$ is a sufficient enlargement for $l_2^n$. Since
$A$ is a sufficient enlargement of minimal volume then $M=A$. 
Let $N=A\cap$lin$\{w_1, w_3,w_4,\dots,w_n\}$.
It is easy to see that $N$ is a sufficient enlargement
for $l_2^{n-1}$ and vol$_nA=2$vol$_{n-1}N$. Hence
vol$_{n-1}N=2^{n-1}$. By induction hypothesis $N$
is a cube. Hence $A$ is also a cube.
$\square$
\enddemo


\Refs
\widestnumber\key{W}

\ref\key 1
\by  Y.Gordon, M.Meyer and A.Pajor
\paper Ratios of volumes and
factorization through $l_\infty$
\jour Illinois J. Math.
\vol 40
\yr 1996
\pages 91--107
\endref

\ref\key 2 
\by B.Gr\"unbaum
\paper Projection constants
\jour Trans. Amer. Math. Soc. 
\vol 95
\yr 1960
\pages 451--465
\endref


\ref\key 3
\by M.I.Kadets and M.G.Snobar
\paper Certain functionals on the Minkowski compactum
\jour Math. Notes
\vol 10
\yr 1971
\pages 694--696
\endref
 
\ref\key 4
\by J.Lindenstrauss and L.Tzafriri
\book Classical Banach Spaces
\bookinfo Lecture Notes in Math., v.338
\publ Springer--Verlag
\yr 1973
\publaddr Berlin
\endref


\ref\key 5
\by M.I.Ostrovskii
\paper Generalization of projection constants:
sufficient enlargements
\jour Extracta Math.
\yr 1996
\vol 11
\pages 466--474
\endref

\ref\key 6
\by R.Schneider
\book Convex Bodies: the Brunn--Minkowski Theory
\bookinfo Encyclopedia of Mathematics and its Applications
\vol 44
\publ Cambridge University Press
\publaddr Cambridge
\yr 1993
\endref

\ref\key 7
\by N.Tomczak-Jaegermann
\book Banach--Mazur distances and finite-dimensional operator
ide\-als
\publ Longman Scientific\&Technical
\publaddr New York
\yr 1989
\bookinfo Pitman Monographs and Surveys in Pure and
Applied Mathematics
\vol 38
\endref

\endRefs
\enddocument